\documentclass[11pt]{article}
\usepackage{color}
\usepackage{amsmath}
\usepackage{mathrsfs}
\usepackage{amssymb}
\usepackage{amsthm}
\usepackage{epstopdf}
\usepackage{graphicx}
\usepackage{mathtools}
\usepackage[usenames,dvipsnames, table]{xcolor}
\usepackage[hang]{caption}
\usepackage{hyperref}
\usepackage{cite}

\usepackage{oldgerm}  
\usepackage[yyyymmdd,hhmmss]{datetime}

\usepackage{mathabx} 
\def\bee{\begin{enumerate}}\def\eee{\end{enumerate}}
\def\bei{\begin{itemize}}\def\eei{\end{itemize}}
\newcommand{\nco}{\newcommand}
\def\R{\mathbb{R}}
\nco{\red}{\color{red}}
\nco{\blue}{\color{blue}}
\nco{\brown}{\color{Magenta}}

\nco{\magenta}{\color{magenta}}

\nco{\violet}{\color{violet}}
\nco{\olive}{\color{Emerald}}
\nco{\orange}{\color{orange}}
\nco{\redend}{\normalcolor}
\nco{\blueend}{\normalcolor}
\def\inv#1{\frac{1}{#1}}
\def\tr{{\rm tr}\,}

\def\ommit#1{{}}
\def\({\left(}
\def\){\right)}
\def\ie{{\it i.e.,\/}\ }
\def\ie{{\rm i.e.,\/}\ }

\definecolor{cb}{rgb}{.8,.5,0}

\nco{\rnc}{\renewcommand}
\rnc{\title}[1]{{\Large\bf\mbox{}\\\medskip#1\bigskip\medskip\\}}
\rnc{\author}[1]{{\large #1\smallskip\\}}
\nco{\address}[1]{{\em #1\medskip\\}}

\def\diag{{\rm diag \,}}
\def\ii{\mathrm{i\,}}
\nco{\bun}{{\bf 1}}
\def\be{\begin{equation}}\def\ee{\end{equation}}
\def\bea{\begin{eqnarray}}\def\eea{\end{eqnarray}}
\def\bee{\begin{enumerate}}\def\eee{\end{enumerate}}
\def\bei{\begin{itemize}}\def\eei{\end{itemize}}
\def\oh{\frac{1}{2}}
\def\ommit#1{{}}

\def\SU{{\rm SU}}\def\U{{\rm U}}
\def\inv#1{\frac{1}{#1}}
\def\mult{\mathrm{br}}
\def\tr{{\rm tr\, }}

\def\eq=#1{\buildrel #1 \over{=}}

\def\CH{{\mathcal H}}      \def\CK{{\mathcal K}}  \def\CO{{\mathcal O}}
\def\CS{{\mathcal S}}

\def\Ga{\pmb{\alpha}}
\def\Gb{\pmb{\beta}}

\def\bbeta{\bar\beta}
\def\bbeta{\gamma}

\def\diag{{\rm diag \,}}
\def\ii{\mathrm{i\,}}

\def\Z{\mathbb{Z}}

\def\R{\mathbb{R}}

\newtheorem{theorem}{Theorem}

\newtheorem{corollary}{Corollary}

\def\comment#1{{\blue #1}}
\def\comment#1{}
\def\adots{\mathinner{\mkern2mu\raise1pt\hbox{.}\mkern3mu\raise4pt\hbox{.}\mkern1mu\raise7pt\hbox{.}}}  
\def\dddots{\mathinner{\mkern2mu\raise9pt\hbox{.}\mkern3mu\raise3pt\hbox{.}\mkern2mu\raise-4pt\hbox{.}}}

\begin{document}
\begin{titlepage}
\begin{center}
\title{On the Minor Problem\\[8pt]
and Branching Coefficients}
\medskip
\author{Jean-Bernard Zuber}
\address{Sorbonne Universit\'e,  UMR 7589, LPTHE, F-75005,  Paris, France\\ \& CNRS, UMR 7589, LPTHE, F-75005, Paris, France}
\bigskip\bigskip
\begin{abstract}
The Minor problem, namely the study of the spectrum of a principal submatrix of a Hermitian matrix taken 
at random on its orbit under conjugation, is revisited, with emphasis on the use of orbital
integrals and on the connection with branching coefficients in the decomposition of an irreducible
representation of $\U(n)$, resp. $\SU(n)$,  into irreps of $\U(n-1)$, resp. $\SU(n-1)$. 
\end{abstract}
\end{center}

 \end{titlepage}

 \section{Introduction}
What we call the Minor problem deals with the following question:  given an $n$-by-$n$ Hermitian matrix
of given spectrum, what can be said about the eigenvalues of one of its $(n-1)\times (n-1)$ 
principal submatrices?
This question has been thoroughly studied and answered by many authors \cite{Bar,Ner,Olsh,Far}.
As several other such questions, this problem of classical  linear algebra has a counterpart
in the realm of representation theory~\cite{He82,GLS, DuVe}, namely the determination of {\it branching coefficients}
of an irreducible representation (irrep) of $\U(n)$, resp. $\SU(n)$, into irreps of $\U(n-1)$,
resp. $\SU(n-1)$. The aim of this note is to review these questions and to make explicit the link,
by use of orbital integrals. It is thus in the same vein as recent works on the Horn
\cite{Z1, CZ1, CMSZ1, CMSZ2, MS} or Schur-Horn \cite{CZ3} problems.

 This paper is organized as follows. 
In sect. 2, I review the classical Minor problem and recall how it may be rephrased in terms
of $\U(n)$ orbital integrals. This suggests a modification, that will be turn out to be 
natural for the case of $\SU(n)$. 
Sect. 3 is devoted to the issue of branching coefficients for the embeddings
$\U(n-1) \subset \U(n)$ and $\SU(n-1)\subset \SU(n)$. While the former is treated by means 
of Gelfand--Tsetlin triangles and does not give rise to multiplicities, as well known since Weyl~\cite{We}, 
the latter requires a new
technique.  This is where the modified integral introduced in sect. 2 proves useful and is
shown to provide an expression of branching coefficients,  see Theorem \ref{thm2}, which is the main result
of this paper. The use of that formula 
as for the behaviour of branching coefficients under {\it stretching}, \ie dilatation of the 
weights, is briefly discussed in the last subsection.

\section{The classical problem}
\subsection{Notations and classical results}
Let us fix notations: If  $A$ is an $n\times n$ Hermitian matrix with known eigenvalues $\alpha_1 \ge \hdots \ge \alpha_n$, what can be said about the eigenvalues $\beta_1\ge \beta_2\ge \hdots \ge\beta_{n-1}$
of one of its principal $(n-1)\times (n-1)$ minor submatrix (``minor" in short\footnote{In the literature, the word ``minor" refers either to the submatrix or to its determinant. We use here in the former sense.})?\\
A first trivial observation is that if we are interested in the statistics of the $\beta$'s as $A$ is taken 
randomly on its $\U(n)$ orbit $\CO_\alpha$, the choice of the minor among the $n$ possible ones is immaterial, since a permutation of rows and columns of $A$ gives another matrix of the orbit.
\\
A second, less trivial, observation is that the $\beta$'s are constrained by the celebrated
Cauchy--Rayleigh interlacing Theorem:
\be\label{support} \alpha_1 \ge\beta_1\ge \alpha_2\ge \beta_2\ge \hdots \beta_{n-1} \ge \alpha_n \,.\ee
For proofs, see for example \cite{Interlacing,Far}

If $A$ is chosen at random on  its orbit $\CO_\alpha$,  and uniformly in the sense of
the $\U(n)$ Haar measure, what is the probability distribution (PDF) of the $\beta$'s ? 
This question has been answered by  Baryshnikov~\cite{Bar}, see also \cite{Olsh,Far}.
We first observe that the problem is invariant under a global shift of all $\alpha$'s and all
$\beta$'s by a same constant: indeed a translation of $A$ by $a \,\Bbb{I}_n$ shifts by $a$ all its eigenvalues 
as well as all the eigenvalues of any of its principal minors. \\
Let $\Delta$ denote the Vandermonde determinant: $\Delta_n(\alpha)=\prod_{1\le i<j\le n} (\alpha_i-\alpha_j)$ 
and likewise for $\Delta_{n-1}(\beta)$.

\begin{theorem}[Baryshnikov\cite{Bar}] \label{thm1}The PDF of the $\beta$'s on its support (\ref{support}) is given by
\be\label{PDF}P(\beta \,|\, \alpha)= (n-1)! \frac{\Delta_{n-1}(\beta)}{\Delta_n(\alpha)}\,.\ee
\end{theorem}

This result may also be recovered in terms of orbital integrals. 
Let
\be\label{HCn} \CH^{(n)}_\alpha(X)= \int_{\U(n)} dU 
e^{\tr U \alpha U^\dagger X}\ee
where $X\in H_n$, the space of $n\times n$ Hermitian matrices, $dU$ is the normalized Haar measure 
on $\U(n)$, and $\alpha$ stands here for  the diagonal matrix $\diag(\alpha_i)$. 
In terms of the eigenvalues $x_i$ of $X$, we write this orbital integral as $\CH^{(n)}(\alpha; x)$.
 
The orbit $\CO_\alpha$ carries a unique probabilistic measure, the {\it orbital measure}
$\mu_\alpha(dA)$, $A\in \CO_\alpha$, 
whose Fourier transform  is  $\CH^{(n)}_\alpha(\ii X)$
$$ \Bbb{E}(e^{\ii \tr A X})  = \int_{\CO_\alpha} e^{\ii \tr A X} \mu_\alpha(dA)=
\int_{\U(n)} dU e^{\ii \tr U \alpha U^\dagger X} = \CH^{(n)}_\alpha(\ii X)\,.$$

Let $\Pi$ be the projector of $H_n$ into $H_{n-1}$ that maps $A\in H_n$ onto 
its upper $(n-1)\times (n-1)$ minor submatrix $B$. 
{According to} the observation that  the Fourier transform of the projection {of the orbital measure} is the restriction 
 of the Fourier transform~\cite {Far}, the characteristic function of $B$ is
 $\phi_B(Y)= \phi_A(X_0)$, with $ X_0=\Pi(X)=\begin{pmatrix} Y&0\\0&0 \end{pmatrix}\in H_n$,
 $Y\in H_{n-1}$,  from which the PDF of $B$ is obtained by inverse Fourier transform 
 $$P(B|A)=\inv{(2\pi)^{(n-1)^2}}  \int_{H_{n-1}} dY e^{-\ii \tr Y B} \int_{\U(n)}  dU e^{\ii \tr  U\alpha U^\dagger X_0}\,.$$
  After reduction to eigenvalues,\footnote{The factor $(n-1)!$ comes from the fact that we are restricting the $\beta$'s 
  to 
  the dominant sector 
 $\beta_1\ge \beta_2\ge \cdots \ge \beta_{n-1}$\,.}
\comment{$$P(\beta \,|\, \alpha)=\inv{(2\pi)^{(n-1)^2}}\frac{(2\pi)^{(n-1)(n-2)}}{(\prod_{1}^{n-1}p!)^2}
{ (n-1)!}
\ \Delta_{n-1}^2(\beta) 
\int_{\R^{n-1}}dx\, \Delta_{n-1}^2(x) 
 \CH^{(n)}(\alpha ;\ii (x,0))  \CH^{(n-1)}(\beta; \ii x)^* \,.$$}
\be P(\beta \,|\, \alpha)=\frac{(n-1)!}{(2\pi)^{n-1} (\prod_{1}^{n-1}p!)^2}
\ \Delta_{n-1}^2(\beta) \int_{\R^{n-1}}dx\, \Delta_{n-1}^2(x) 
 \CH^{(n)}(\alpha ;\ii (x,0))  \CH^{(n-1)}(\beta; \ii x)^* \,.
\ee
(In physicist's parlance, this is the {\it overlap} of the two orbital integrals.)
Here and below, for $x\in \R^{n-1}$, $(x,0)$ denotes the corresponding vector in $\R^n$. 

Making use of the explicit expressions known for $\CH^{(n)}(\alpha;x)$~\cite{HC,IZ}, we find
\bea \label{Pba}
\!\!\!\!\!\!\!\!\!P(\beta \,|\, \alpha)\!\!\!&=&\!\!\! C 
\ \Delta_{n-1}^2(\beta) \int_{\R^{n-1}} dx \Delta_{n-1}^2(x) 
\frac{\det \left(e^{\ii \alpha_i (x,0)_j}\right)_{i,j=1,\cdots,n}}{\Delta_{n}(\alpha)\Delta_{n}((x,0))}  
\frac{\det \left(e^{-\ii \beta_i x_j}\right)_{i,j=1,\cdots,n-1}}{\Delta_{n-1}(\beta)\Delta_{n-1}(x)} \\
\!\!\!&=& \!\!\!
\label{Pbetalf}
\inv{ (2\pi \ii)^{n-1}} \ \frac{\Delta_{n-1}(\beta)}{\Delta_n(\alpha)} \int_{\R^{n-1}} \frac{d^{n-1}x}{x_1 x_2\cdots x_{n-1}}
\det e^{\ii \alpha_i (x,0)_j}\, \det e^{-\ii \beta_i x_j}
\eea
{since the prefactor reads  $$ 
C=  \frac{{(n-1)!} \prod_1^{n-1}p!  \, \prod_1^{n-2}p! \ \ii^{-n(n-1)/2 +(n-1)(n-2)/2}}{(2\pi)^{n-1} (\prod_{1}^{n-1}p!)^2} =\inv{ (2\pi \ii)^{n-1}} $$}
and since $\Delta_{n}((x,0)) =( \prod_{i=1}^{n-1} x_i )\Delta_{n-1}(x)$.
Let's write $P(\beta \,|\, \alpha)=\frac{\Delta_{n-1}(\beta)}{\Delta_n(\alpha)} \CK(\alpha;\beta)$, hence
 \bea \nonumber \label{defKa0} 
   \CK(\alpha;\beta) &=&   
 \frac{(n-1)!}{(2\pi)^{n-1} (\prod_{1}^{n-1}p!)^2}
 \ \Delta_n(\alpha) \Delta_{n-1}(\beta) 
    \int_{\R^{n-1}}dx \Delta_{n-1}^2(x) 
 \CH^{(n)}(\alpha;\ii (x,0))  \CH^{(n-1)}(\beta;\ii x)^*\\
\label{defKa}  &=& 
\inv{ (2\pi \ii)^{n-1}}  \int_{\R^{n-1}} \frac{d^{n-1}x}{x_1 x_2\cdots x_{n-1}}
\det \big(e^{\ii \alpha_i (x,0)_j}\big)_{\scriptscriptstyle 1\le i,j \le n}\, \det \big(e^{-\ii \beta_i x_j}\big)_{\scriptscriptstyle 1\le i,j\le n-1}
 \eea 
 in analogy with the introduction of the ``volume functions" in the Horn and Schur problems~\cite{CMSZ1,CMSZ2,CZ3}. 
The function  $\CK(\alpha;\beta)$ is then, as in these similar cases, a linear 
combination of products of Dirichlet integrals: $\mathrm{PP}\int \frac{e^{\ii \alpha t}}{t}=\ii \pi \varepsilon(\alpha)$,
with $\varepsilon$ the sign function. Thus $\CK(\alpha;\beta)$ must be \comment{[a piecewise 
polynomial function of degree 0, hence]} a piecewise constant function,
supported by the product of intervals given by the interlacing theorem (\ref{support}).

By making use of the integral form of the Binet--Cauchy formula, (see \cite{Far}),
namely
$$ \int_{\R^k} \det ( f_i(t_j))_{\scriptscriptstyle{1\le i,j\le k}} \ \det( g_i(t_j))_{\scriptscriptstyle 1\le i,j\le k} = k! \det \(\int_R f_i(t) g_j(t)\)_{\scriptscriptstyle 1\le i,j\le k}\,,$$
with here $k=n-1$, $f_i(t) = \frac{1}{t} (e^{\ii \alpha_i t}-e^{\ii \alpha_n t})$, $g_i(t)=e^{-\ii \beta_i t}$, 
we find
\bea\label{newK}  
  \CK(\alpha;\beta) &=& \frac{(n-1)! }{ (2\pi \ii)^{n-1} } 
   \det \int_\R \frac{dt}{t} ( e^{\ii (\alpha_i-\beta_j) t}-e^{\ii (\alpha_n-\beta_j) t} ) \\  \label{altK1}
&=&\frac{(n-1)! }{ 2^{n-1}}\, 
\det \big(\varepsilon(\alpha_i-\beta_j)-
\varepsilon(\alpha_n-\beta_j) \big)_{1\le i,j\le n-1}\\
\label{altK}   &=& 
\frac{(n-1)!}{2^{n-1}} \det \begin{pmatrix}
\varepsilon(\alpha_1-\beta_1) & \cdots & \varepsilon(\alpha_1-\beta_{n-1}) & 1\\
\varepsilon(\alpha_2-\beta_1) & \cdots & \varepsilon(\alpha_{2}-\beta_{n-1}) & 1\\
\vdots && &\vdots \\
\varepsilon(\alpha_n-\beta_1) & \cdots & \varepsilon(\alpha_{n}-\beta_{n-1}) & 1   
\end{pmatrix}\,.\eea

Equ.  (\ref{altK1}) just reproduces a result by Olshanksi \cite{Olsh},
since the difference $\big(\varepsilon(\alpha_i-\beta_j)-
\varepsilon(\alpha_n-\beta_j) \big)$ appearing there is nothing else than twice the
characteristic function of the interval $[\alpha_n, \alpha_i]$, denoted 
$M_2(\beta_j; \alpha_n,\alpha_i)$ in \cite{Far}.
Finally, it may be shown that
the determinant in (\ref{altK}) equals $2^{n-1}$ times the characteristic function of (\ref{support}), 
so that the piecewise constant function $\CK$  is just $(n-1)!$ on its support, in agreement with (\ref{PDF}),
see \cite{Olsh, Far}.

{\bf Remark}. The previous considerations extend to projections of matrix $A$ onto a
smaller minor $k\times k$ submatrix, see \cite{Ner, Olsh, Far}.
\def\CKp{\widebar \CK}

\subsection{A modified integral} 
\label{modpb}
In this section, we introduce a modification of the integral (\ref{defKa}) that will be suited in our later
study of $\SU(n)$ branching rules.
We first change variables in $\CK$, introducing the spacings $\bbeta_i=\beta_i-\beta_{n-1}$,
$i=1,\cdots,n-2$  and rewrite (\ref{defKa}) as
$$  \CK(\alpha;\beta)=  \inv{ (2\pi \ii)^{n-1}  }  \int_{\R^{n-1}} \frac{d^{n-1}x}{x_1 x_2\cdots x_{n-1}}  e^{-\ii \beta_{n-1} \sum_{j=1}^{n-1}x_j}\,
\det \big( e^{\ii \alpha_i (x,0)_j}\big)_{\scriptscriptstyle1\le i,j\le n}\, \det \big(e^{-\ii \bbeta_i  x_j}\big)_{\scriptscriptstyle
1\le i,j\le n-1}\,,
$$
where by convention $\bbeta_{n-1}=0$. We then integrate  
over $\beta_{n-1}$ (while introducing a $1/(n-1)!$ prefactor for later convenience),
and define  
\bea  \label{defKp} \!\!\!\!\!
\CKp(\alpha;\bbeta)\!\!\! &:=&\!\!\!  \inv{(n-1)!} \int d\beta_{n-1} \CK(\alpha;\bbeta+\beta_{n-1})\\
&=&\!\!\!  \nonumber 
\frac{ (2\pi) }{(n-1)! (2\pi \ii)^{n-1} }   
\int_{\R^{n-1}} \frac{d^{n-1}x}{x_1 x_2\cdots x_{n-1}}  \ \delta\Big(\sum_1^{n-1} x_j\Big)
\det \big( e^{\ii \alpha_i (x,0)_j}\big)_{\scriptscriptstyle 1\le i,j\le n}\, \det \big(e^{-\ii \bbeta_i x_j}\big)_{\scriptscriptstyle 1\le i,j\le n-1}
 \,.
 \eea
 {\bf Remark}. Although it depends only on the spacings  $\bbeta$, the function 
 $\CKp$ {\it is not}  directly related to the PDF of the spacings in the original Minor problem. Its
 introduction is rather motivated by its connection with the $\mathfrak{su}(n)$ Lie algebra, see below sect. \ref{A-CK}. 
 
Thus integrating $\CK$ over $\beta_{n-1}$  amounts to considering   a modified integral, where in (\ref{defKa}) we integrate on the 
$(n-2)$-dimensional 
hyperplane $\sum_{i=1}^{n-1} x_i=0$.  Hence an alternative definition of $\CKp$ is
\be\label{later} 
\CKp(\alpha;\bbeta) = 
 \frac{\Delta_n(\alpha) \Delta_{n-1}(\bbeta) }{(2\pi)^{n-2} (\prod_{1}^{n-1}p!)^2}
  \int_{\R^{n-1}}dx \,\delta\Big(\sum_1^{n-1} x_j\Big) \,\Delta_{n-1}^2(x) 
 \CH^{(n)}(\alpha;\ii (x,0))  \CH^{(n-1)}(\bbeta;\ii x)^*\,,
 \ee
an expression that we use later in sect. \ref{A-CK}.
{A more explicit expression is
 \be 
\label{modKb} \!\!\!\!\!\CKp(\alpha;\bbeta) 
={\textstyle \frac{2\pi }{  (2\pi \ii)^{n-1}  (n-1)!} }\  \int_{\sum_1^{n-1} \!\!x_i=0}\ \,
\frac{dx}{x_1 x_2\cdots x_{n-1}}
\det \left({\scriptstyle e^{\ii \alpha_i x_1}, e^{\ii \alpha_i x_2},\cdots, e^{\ii \alpha_i x_{n-1}}, 1}\right)_{\scriptscriptstyle 1\le i \le n-1} \,\det\big( e^{-\ii \bbeta_i  x_j}\big)_{\scriptscriptstyle 1\le i,j\le n-1}
 \ee
and we note that, because of  the constraint $\sum_{i=1}^{n-1} x_i=0$, this expression is 
invariant by a global shift of all $\alpha_i$. We may use that invariance to choose $\alpha_n=0$, a choice that 
will be natural in the application to $\SU(n)$ representations. 
We conclude that $\CKp(\alpha,{\bbeta})$ is a function of two sets of variables, 
a $n$-plet $\alpha$ with $\alpha_n=0$, and a $(n-1)$-plet ${\bbeta}$ with ${\bbeta}_{n-1}=0$.}
{As is clear from (\ref{defKp}), $\CKp$ may be extended to a function of the unordered $\alpha$'s and ${\bbeta}$'s, 
odd under the action of the symmetric  group $\CS_{n}$, \ie the $\SU(n)$ Weyl group, acting on $\alpha$ by
$ w(\alpha)_i=\alpha_{w(i)} - \alpha_{w(n)}$, $i=1,\cdots, n$, 
$w\in \CS_n$, and likewise odd under the action of
$\CS_{n-1}$  on  ${\bbeta}_i$, $i=1,\cdots, n-1$. }

 Expanding the two determinants and using once again the Dirichlet integrals $\mathrm{PP}\int \frac{e^{\ii a t}}{t^r}=\ii \pi \frac{(\ii a)^{r-1}}{(r-1)!}\varepsilon(a)$, 
   one finds that $\CKp(\alpha; {\bbeta})$, a combination of convoluted box splines,  is a piece-wise linear function
 of differentiability class $C^0$.
 Its support is the polytope defined by the inequalities (recall that by convention $\alpha_n=0$)
 \be\label{polytope}\max_{1\le i \le n-1}(\alpha_{i+1}-{\bbeta}_i) \le \min_{1\le i \le n-1}(\alpha_i-{\bbeta}_i)
 \ee
 that guarantee that there exist $\beta_{n-1}$ satisfying  the simultaneous inequalities
 (\ref{support}), \ie
\be \label{condbetan} \alpha_{i+1}\le \beta_i= {\bbeta}_i+\beta_{n-1} \le \alpha_{i}\,,\ 
 \mathrm{for\ all}\ i=1,\cdots, n-1\,. \ee
 The maximal value  (in the dominant sector) of $\CKp$, for fixed $\alpha$, is readily derived from 
 (\ref{defKp}), where we are integrating the function $\CK/(n-1)!$ equal to 1 on its support, 
 over $\beta_{n-1}$, subject to the $n-1$ conditions  (\ref{condbetan}),  hence 
 \be\label{maxCKp} \max_{{\bbeta}} \CKp(\alpha;{\bbeta}) 
 =\min_{1\le i\le n-1} (\alpha_{i}-\alpha_{i+1})\,.\ee
 \\

\begin{figure}[htb]
  \centering
       \includegraphics[width=15pc]{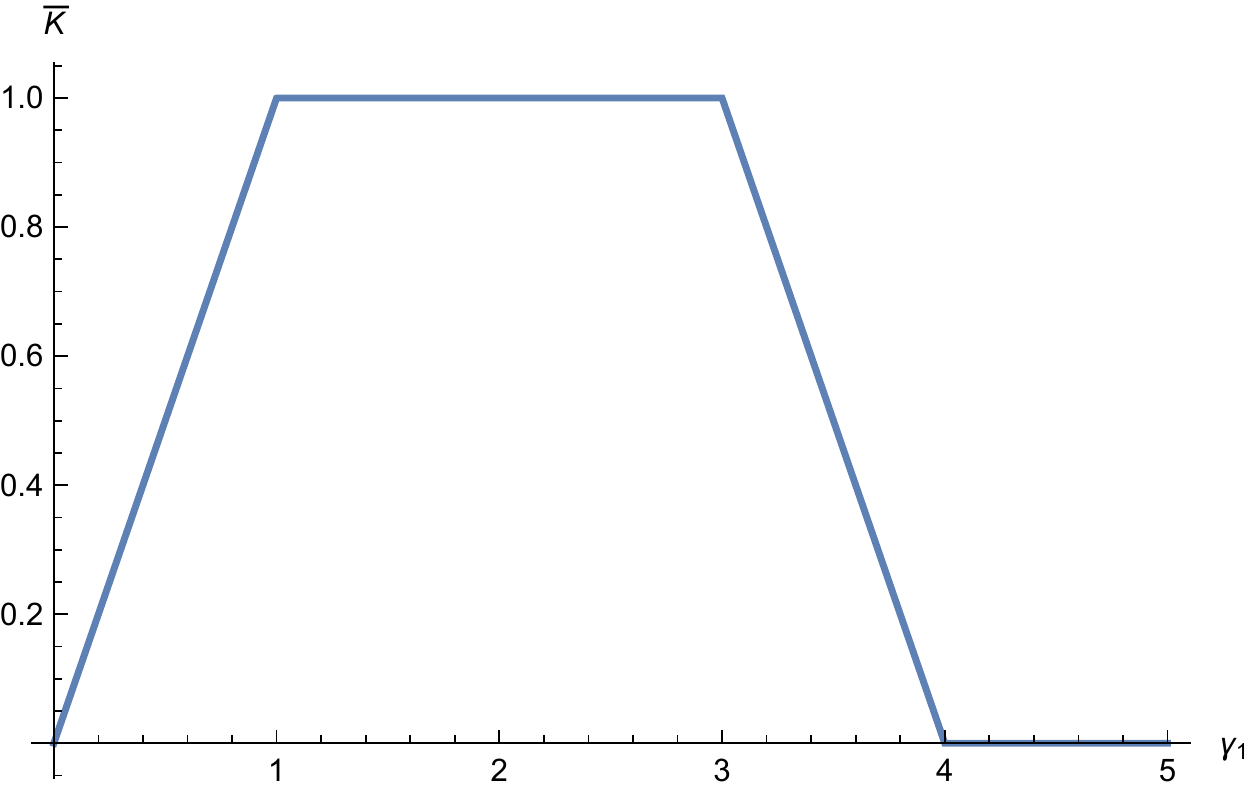}\caption{\label{figKp3} The $\CKp$ function for $n=3$,  $\alpha=\{4,1,0\}$  and ${\bbeta}_1\ge {\bbeta}_2=0$.   }
       \end{figure}

\begin{figure}[htb]
  \centering
       \includegraphics[width=17pc]{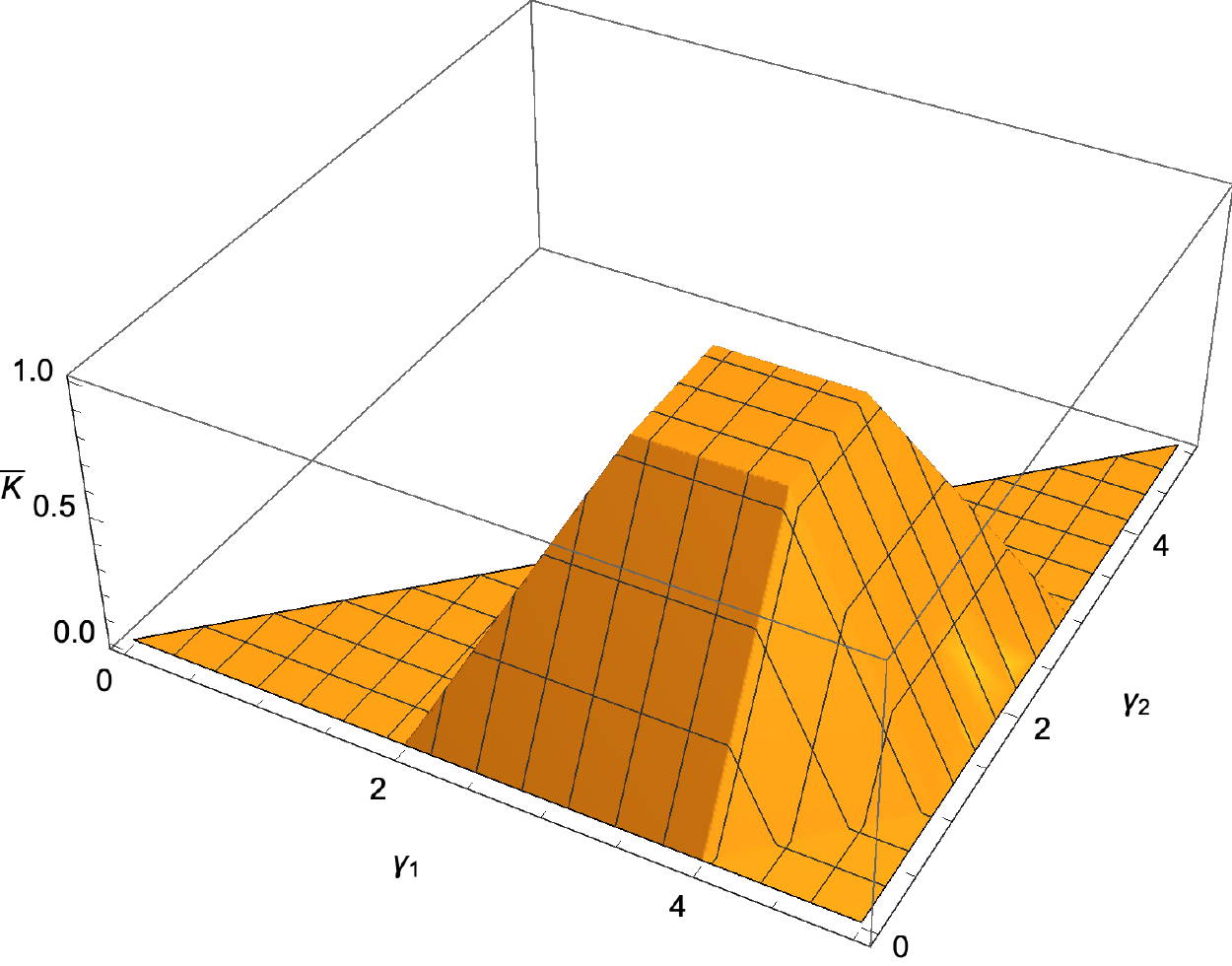}
             \caption{\label{figKp4}
       The $\CKp$ function  for $n=4$  and $\alpha=\{5,3,1,0\}$,
       in the ${\bbeta}_1\ge  {\bbeta}_2\ge {\bbeta}_3=0$ sector.       }
       \end{figure}

 Let $\epsilon_w$ denote the signature of permutation $w\in \CS_n$.
 For  $n=3$ the function $\CKp$ reads 
   \bea \label{Kp3}
\CKp(\alpha;{\bbeta})&=&
    \oh\Big(\left|\alpha _1-{\bbeta} _1\right| 
  -\left|\alpha _1-\alpha _2-{\bbeta} _1\right| 
   -\left|\alpha_2-{\bbeta} _1\right| 
   \Big)
   - ({\bbeta}_1\mapsto -{\bbeta}_1)\\
   \nonumber &=& \oh \sum_{w\in \CS_3} \epsilon_w |w(\alpha)_1-{\bbeta}_1|
   \eea
 which is an odd continuous function of ${\bbeta}_1$, vanishing for ${\bbeta}_1\notin(-\alpha_1,\alpha_1)$,
 constant and equal to its extremum value $\pm \min( \alpha_1-\alpha_2,\alpha_2)$ for $|{\bbeta}_1|\in [\min(\alpha_1-\alpha_2,\alpha_2), \max(\alpha_1-\alpha_2,\alpha_2))$,
 and linear in between, see Fig. \ref{figKp3}. 
  \\
{ For $n=4$, let $ w(\bar\alpha)_i=\alpha_{w(i)} - \alpha_{w(n)}$, 
   \be   \psi(\alpha;{\bbeta}):=  \varepsilon(\alpha_1-{\bbeta}_1) 
   \Big( |\alpha_2-{\bbeta}_2| 
   -|\alpha_3-{\bbeta}_2|  -|\alpha_1-\alpha_2-{\bbeta}_1+{\bbeta}_2|   
   +|\alpha_1-\alpha_3-{\bbeta}_1+{\bbeta}_2|  
   \Big)\ee
   then
    \be \label{CKp4}\CKp(\alpha;{\bbeta}) =\inv{8} \sum_{w\in \CS_4}  \epsilon_w \psi(w(\alpha);{\bbeta}) \ee
$\CKp(\alpha;{\bbeta})$  has  a support in the dominant sector defined by the inequalities
 \be\label{supportKp4}  
   (\alpha_2-\alpha_3)\le {\bbeta}_1\le \alpha_1\,, \qquad 0\le {\bbeta}_{2} \le   \alpha_2 \,,\qquad
0\le  {\bbeta}_1-{\bbeta}_2\le  \alpha_1-\alpha_3\,, \ee
and a maximal value  equal to $\min(\alpha_1-\alpha_2, \alpha_2-\alpha_3, \alpha_3)$. 
Its  graph 
has an Aztec pyramid shape, see Fig. \ref{figKp4}.}

\section{The ``quantum" problem}
In this section, we consider the restriction of the group $\U(n)$, resp. $\SU(n)$, to its subgroup
$\U(n-1)$, resp. $\SU(n-1)$, and the ensuing decomposition of their representations.
For definiteness, the restriction of $\SU(n)$ to $\SU(n-1)$ we have in mind results 
from projecting out  the  simple root $\Ga_{n-1}$ in the dual of the Lie algebra $\mathfrak{su}(n)$, and likewise for $\U(n)$.

 \subsection{Gelfand--Tsetlin patterns} 
 \label{GTpatterns}
Just like in the cases of the Horn or of the Schur problem, the Minor problem is
the classical counterpart of a ``quantum" problem in representation theory.  Given a
highest weight (h.w.) irreducible representation (irrep) $V^{(n)}_\alpha$ of $\U(n)$, which
irreps $V^{(n-1)}_\beta$
of $\U(n-1)$ occur and with which multiplicities, in the restriction of $\U(n)$ to $\U(n-1)$? 
That problem too is well known, is important in physical applications (see for example
\cite{Geo,DFMS}), and may be solved by a variety of methods. Here we first recall how to make  use 
of Gelfand--Tsetlin triangles, \ie triangular patterns

 $$\begin{array}{rlcrlcrlcrlcrl} 
x_1^{(n)}=\alpha_1&\cdots & & \cdots&  &\cdots&x_n^{(n)}=\alpha_n\\
& x_1^{(n-1)}&\cdots & & \cdots &x_{n-1}^{(n-1)}& \\
&& \ddots & & \adots &&& \\
&&& x_1^{(1)}& && 
  \end{array}$$
  subject to  the inequalities
  \be\label{GTineq} 
    x_{i}^{(j+1)} \ge x_{i}^{(j)} \ge  x_{i+1}^{(j+1)} 
    \ ,\qquad {1 \le i,j \le n-1 }\,.\ee

In the present context, the $\alpha$'s denote the lengths of the rows of the Young diagram associated
with the irrep $V^{(n)}_\alpha$. 
The number of solutions of the inequalities (\ref{GTineq}) gives the dimension of the irrep
$$ \dim(V^{(n)}_\alpha) =\#\{ x_i^{(j)} | \mathrm{solutions\ of\ }  (\ref{GTineq}) \}\,.$$
The values $\beta_i=x_{i}^{(n-1)}$, $1\le i\le n-1$, appearing in the second row of the triangle
give the lengths of rows of the Young diagrams of the possible representations $V^{(n-1)}_\beta$ of
$\U(n-1)$.
Given those numbers, the number of solutions $x_i^{(j)}$, $1\le i, j \le n-2$ satisfying 
(\ref{GTineq}) is the dimension of the representation of $\U(n-1)$. Thus we have the
sum-rule
\bea
\dim(V^{(n)}_\alpha) &=&\#\{ x_i^{(j)} | \mathrm{solutions\ of\ }  (\ref{GTineq}),\ x^{(n)}=\alpha  \}\ \\
\nonumber &=&\sum_\beta 
   \#\{ x_i^{(j)} | \mathrm{solutions\ of\ }  (\ref{GTineq}),\, \ x^{(n)}=\alpha ,\, \ x^{(n-1)}=\beta   \}\ = \sum_\beta \dim(V^{(n-1)}_\beta) \,,\eea
   which is consistent    with the multiplicity 1 of each $V^{(n-1)}_\beta$ appearing in the decomposition,
   {a classical result in representation theory~\cite{We,Mu,Zhe}, see also chapter 8 in \cite{G-W}}.
   Thus one sees that the $\beta$'s satisfy the inequalities (\ref{support})
   and one may say that the branching coefficient, equal to 0 or 1, is given
   \be \mult_\alpha(\beta) = K(\alpha; \beta) \ee
   with the convention that the discontinuous function $\CK$ is assigned the value $1$ throughout
     its support, including its boundaries.

   \medskip
   Going from $\U(n)$ to $\SU(n)$, we have to restrict to Young diagrams with less than $n$ rows, 
   or equivalently, to reduce Young diagrams with $n$ rows by deleting all columns of height $n$.
   Starting from an irrep of $\SU(n)$, we apply to it the procedure above, and then remove the columns of height $n-1$.
   
   Example in $\SU(3)$. Take for $\alpha$  the adjoint representation, \ie $\alpha=\{2,1,0\}$.
   The possible $\beta$ satisfying (\ref{support}) are written in red in what follows
   \bea & 2\ge {\red 2} \ge 1  \ge {\red 1} \ge 0 & \quad \ie \beta=\{2,1\}\equiv \{1,0\}\\
           & 2\ge {\red 1} \ge 1  \ge {\red 1} \ge 0 & \quad \ie \beta=\{1,1\}\equiv \{0,0\}\\
           & 2\ge {\red 2} \ge 1  \ge {\red 0} \ge 0 & \quad \ie \beta=\{2,0\}   \\
           & 2\ge {\red 1} \ge 1  \ge {\red 0} \ge 0 & \quad \ie \beta=\{1,0\}
           \eea
   where two $\beta$ are regarded as equivalent if their Young diagrams differ by a number of
   columns of height     $n-1=2$. Hence in $\SU(2)\subset \SU(3)$, we write
      \be\label{example3}\mult_\alpha(\beta)=   1,2,1 \quad \mathrm{for}\ \beta=  \{0,0\},\ \{1,0\},\ \{2,0\}\ee
       and we check the sum-rule on dimensions: $8=1 + 2\times 2 + 3$. 
       { Note that removing columns of height $n-1$ in the Young diagram associated with $\beta$
       amounts to focusing on spacings between 
       the $\beta$'s, which points to the relevance of our function $\CKp$.}
   
   \bigskip
   To summarize, in the ``quantum", $\U(n)$-representation theoretic, problem, the $\beta$'s are
   the integer points interlacing the $\alpha$'s and come with multiplicity 1, 
   while in the $\SU(n)$ case, non trivial multiplicities may occur and the interlacing property 
   no longer applies. In the next subsection, we show that the latter multiplicities are 
   given by the function $\CKp$ of sect. \ref{modpb}.

\subsection{A $\CKp-\mult$ relation}
\label{A-CK}

The multiplicities occurring in the $\SU(n-1)\subset \SU(n)$ problem may be expressed in terms of characters
by the integral 
\be \mult_\alpha({\bbeta})=\int  Dt \, 
\chi^{(n)}_\alpha(e^{\ii(t,0)}) \big(\chi^{(n-1)}_{\bbeta}(e^{\ii t})\big)^*\ee
which computes the projection of the  $\SU(n)$ character $\chi^{(n)}_\alpha$ restricted to the Cartan torus of 
$\SU(n-1)$ onto the $\SU(n-1)$ character  $\chi^{(n-1)}_{\bbeta}$. 
There, $Dt$ stands for the Haar measure on the 
Cartan torus $\Bbb{T}_{n-1}$ of $\SU(n-1)$
$$ Dt=\frac{|\widehat{\Delta}_{n-1}(e^{\ii t})|^2}{(2\pi)^{n-2} (n-1)! } dt 
$$
where we use the notations 
$$ 
\widehat{\Delta}_{n-1}(e^{\ii t}):= \prod_{\Ga >0} \big(e^{\ii \langle \Ga,t \rangle/2}-e^{-\ii \langle \Ga,t \rangle/2}\big) \,, \qquad\mathrm{and}\qquad
\Delta_{n-1}(t):= \prod_{\Ga >0}   \langle\Ga,t\rangle\,,$$
$\Ga$ the positive roots of $\mathfrak{su}(n-1)$, and $dt$ is the Lebesgue measure on 
$\Bbb{T}_{n-1}$. 

\medskip
\bigskip
{\begin{theorem}\label{thm2} The branching coefficient, that gives the multiplicity of  
the  irrep of $\SU(n-1)$ of h.w. ${\bbeta}$ 
in the decomposition  of
the irrep of $\SU(n)$ of h.w. $\alpha$, is 
\be \label{K-br}  \mult_\alpha({\bbeta})=\CKp(\alpha+\rho_n;{\bbeta}+\rho_{n-1})\ee
with $\rho_n$ the Weyl vector of the algebra  $\mathfrak{su}(n)$, and $\rho_{n-1}$ that of $\mathfrak{su}(n-1)$.
\end{theorem}}

\noindent {\it Proof}. 
We recall Kirillov's relation between a $\SU(n)$ character and the orbital integral:
\be\label{kri}   
 \chi^{(n)}_\alpha(e^{\ii t})= \dim V_\alpha  \frac{\Delta_n(\ii t)}{\widehat{\Delta}_n(e^{\ii t})} 
 \CH^{(n)}(   {\alpha+\rho_n}\,; \ii t) 
\ee
with 
$ \dim V_\alpha =\frac{\Delta_n(\alpha+\rho_n)}{\Delta_n(\rho)}$.
Plugging   in (\ref{later}) the expression (\ref{kri}) and the analogous one for $\mathfrak{su}(n-1)$
 leads to
\bea \nonumber
\!\!\!\!\!\CKp(\alpha+\rho_n;{\bbeta}+\rho_{n-1}) 
\!\!\! &=&\!\!\!    \frac{ \Delta_n(\alpha+\rho_n)  \Delta_{n-1}({\bbeta}+\rho_{n-1}) }{(2\pi)^{n-2} (\prod_{1}^{n-1}p!)^2} 
\int_{\R^{n-2}} \!\!\!  \!\!\! dt \,\Delta_{n-1}^2(t) 
 \CH^{(n)}( {\alpha+\rho_n} ; \ii(t,0))  \CH^{(n-1)}( {{\bbeta}+\rho_{n-1}};\ii t)^*
  \\
 \nonumber \!\!\!\!\!\ &=&
 \frac{  \prod_{p=1}^{n-1}p!  \prod_{p=1}^{n-2}p!  }{(2\pi)^{n-2} (\prod_{1}^{n-1}p!)^2} 
\int_{\R^{n-2}} \!\!\!  dt\,   \Delta_{n-1}^2(t) 
 \frac{\widehat{\Delta}_n(e^{\ii(t,0)})  \widehat{\Delta}_{n-1}(e^{\ii t})^*}{\Delta_n(\ii (t,0))  \Delta_{n-1}( \ii t)^*}
  \chi^{(n)}_\alpha(e^{\ii (t,0)})   \chi^{(n-1)}_{\bbeta}(e^{\ii t})^*
  \\
   \nonumber &=&  \frac{\ii^{-(n-1)}}    {(2\pi)^{n-2} (n-1)!}  
 \int_{\R^{n-2}}     \!\!\!  dt \,  |\widehat{\Delta}_{n-1}(e^{\ii t})|^2
    \chi^{(n)}_\alpha(e^{\ii (t,0)})   \chi^{(n-1)}_{\bbeta}(e^{\ii t})^* \frac{\Delta_{n-1}(t)}{\Delta_n((t,0))}
    \frac{\widehat{\Delta}_n(e^{\ii(t,0)})}{ \widehat{\Delta}_{n-1}(e^{\ii t})}
      \\
 \nonumber &=&   \ii^{-(n-1)} \,  \int_{\Bbb{T}_{n-1}}  \!\!\!  Dt\, 
        \sum_{\delta\in 2\pi Q^\vee} 
     \chi^{(n)}_\alpha(e^{\ii (t+\delta,0)})   \chi^{(n-1)}_{\bbeta}(e^{\ii (t+\delta)})^* 
         \frac{\Delta_{n-1}(t+\delta)}{\Delta_n((t+\delta,0))}
    \frac{\widehat{\Delta}_n(e^{\ii(t+\delta,0)})}{ \widehat{\Delta}_{n-1}(e^{\ii (t+\delta)})}
\eea
where the integration 
 is now carried out on the Cartan torus of $\SU(n-1)$,  $\Bbb{T}_{n-1}=\R^{n-2}/(2\pi Q^\vee)$,
$Q^\vee$ the $(n-2)$-dimensional coroot lattice of  $\SU(n-1)$, which is, in the present simply laced case, 
isomorphic to the root lattice. Only the ratio $ \frac{\Delta_{n-1}(t+\delta)}{\Delta_n((t+\delta,0))}$
depends on $\delta$  and the summation can be carried out with the result that 
\be\label{niceidentity} 
\ii^{-(n-1)} \frac{\widehat{\Delta}_n(e^{\ii(t,0)})}{ \widehat{\Delta}_{n-1}(e^{\ii t})} \sum_{\delta\in 2\pi Q^\vee}   \frac{\Delta_{n-1}(t+\delta)}{\Delta_n((t+\delta,0))}=1 \,.\ee
Indeed, if we write $(t,0)$ in the $\mathfrak{su}(n)$ root basis: $(t,0)=\sum_{j=1}^{n-2} a_j \Ga_j$ 
(with no component on $\Ga_{n-1})$, 
\be\label{telprod} \frac{\widehat{\Delta}_n(e^{\ii(t,0)})}{ \widehat{\Delta}_{n-1}(e^{\ii t})}  = - (2\ii)^{n-1} \sin \frac{a_1}{2} 
\(\prod_{i=1}^{n-3} \sin\frac{a_{i+1}-a_i}{2} \)\, \sin \frac{a_{n-2}}{2} \,,\ee
(on which it is clear that it is invariant under $a_i\mapsto a_i +p_i(2\pi)$, $p_i\in \Z$), 
while 
\be \label{telprod} \frac{\Delta_n((t,0))}{\Delta_{n-1}(t)}= 
-{a_1} \(\prod_{i=1}^{n-3} (a_{i+1}-a_i)\) \, {a_{n-2}} \ee
and the identity (\ref{niceidentity}) follows from a repeated use of
$$ \sum_{p=-\infty}^\infty \inv{(a+2\pi p)(b-a-2\pi p)}
=
\inv{2b }\, \frac{\sin\frac{b}{2}}{\sin\frac{a}{2}  \sin\frac{b-a}{2}}\,$$
in the telescopic product (\ref{telprod}).
\comment{[To wit, 
\bea &&\sum_{p_1,p_2} \inv{(a_1+2\pi p_1) (a_2+2\pi p_2-a_1-2\pi p_1)(a_3-a_2-2\pi p_2)}\\
&=&\oh \frac{\sin\frac{a_2}{2}}{\sin\frac{a_1}{2}  \sin\frac{a_2-a_1}{2}}\sum_{p_2}\inv{(a_2+2\pi p_2)(a_3-a_2-2\pi p_2)}\\
\nonumber  &=& \inv{2^2\, a_3}\,\, \frac{\sin\frac{a_2}{2}  \sin\frac{a_3}{2}}{\sin\frac{a_1}{2}  
\sin\frac{a_2}{2} \sin\frac{a_2-a_1}{2} \sin\frac{a_3-a_2}{2}}\eea
etc.  At the last step :
$$\frac{\sin\frac{a_{n-2}}{2}}{\sin\frac{a_{n-2}-a_{n-3}}{2} }\sum_{p_{n-2}} \oh \inv{(a_{n-2}+2\pi p_{n-2})^2}
=\frac{\sin\frac{a_{n-2}}{2}}{\sin\frac{a_{n-2}-a_{n-3}}{2} } \inv{8 \sin^2\frac{a_{n-2}}{2}}
=\inv{8   \sin\frac{a_{n-2}-a_{n-3}}{2} \sin\frac{a_{n-2}}{2}}
$$
Final coefficient in (\ref{niceidentity} ): $ \ii^{-(n-1)} (2\ii)^{n-1}  \inv{2^{n-3+2}}=  1$, \quad phew!]
}
\hfill $\Box$
\\[8pt]

{\bf Remark}. The proof above follows closely similar proofs in \cite{CMSZ1,CZ3} that relate the classical
Horn or Horn--Schur problems to the computation of Littlewood--Richardson or Kostka coefficients. However,
in contrast with those cases, here the r.h.s is a single term, 
 rather than a linear combination  involving a convolution. 
\\[8pt]

{Together with the results of the end of sect. \ref{modpb}, Theorem  \ref{thm2} has immediate consequences:
\begin{corollary} The number of irreps of $\SU(n-1)$ appearing in the decomposition of the irrep of $\SU(n)$ of h.w. $\alpha$ 
(with $\alpha_n=0$) is equal to the number of integer points in the polytope defined by the inequalities (\ref{polytope}), where $ \alpha$ is 
changed into $\alpha+\rho_n$, namely
$$ \# \{ {\bbeta} \in \Z_+^{n-1}\ |\  {\bbeta}_1\ge \cdots\ge {\bbeta}_{n-1}\ge 0\ ,\quad  \max_{1\le i \le n-1}(\alpha_{i+1}+n-i-1 -{\bbeta}_i) \le \min_{1\le i \le n-1}(\alpha_i+n-i-{\bbeta}_i) \}\,.$$
\end{corollary}}
\noindent Here as before, the $\alpha_i$ are the Young coordinates of the h.w. $\alpha$, \ie, the lengths of the rows of its
Young diagram.

On the other hand,  eq. (\ref{K-br}) together with (\ref{maxCKp}) gives the maximal value of a branching coefficient
of a given $\alpha$
 \be\label{boundbr} \max_{{\bbeta}} \mult_\alpha({\bbeta})= \min_i((\alpha+\rho_n)_i-(\alpha+\rho_n)_{i+1}) =\min_i(\alpha_i-\alpha_{i+1})+1
\ee
and note that $\alpha_i-\alpha_{i+1}$ is just the $i$-th Dynkin component\footnote{Recall  that the Dynkin components of  a
weight are its components in the fundamental weight basis. Hereafter, they are denoted by round brackets.} of the weight $\alpha$.
\begin{corollary} The largest multiplicity (branching coefficient) that occurs in the branching of
an irrep of $\SU(n)$ of h.w. $\alpha$ into irreps of $\SU(n-1)$ is 1 plus the smallest Dynkin 
component of $\alpha$. \end{corollary}

\noindent {{\bf Examples}. 
Take $n=3$ and the example considered in sect. \ref{GTpatterns}.
$\alpha=\{2,1,0\}$, \ie $(1,1)$ in Dynkin components, $\alpha+\rho_3=\{4,2,0\}$, ${\bbeta}\in\{ \{0,0\},\, \{1,0\},\,\{2,0\}\},\ 
{\bbeta}+\rho_2 \in\{ \{1,0\},\, \{2,0\},\, \{3,0\}\}$, 
one finds with the formula (\ref{Kp3}): $\CKp(\alpha+\rho_3;{\bbeta}+\rho_2)= 1,2,1$ in agreement with (\ref{example3}). \\
For $n=4$, take $\alpha=\{6,4,3,0\}$, {\ie $(2,1,3)$ in Dynkin components},
one finds the following decomposition into 18 $\SU(3)$ weights
\bea \{6,4,3,0\} \! \! \! \! \! &=&\! \! \! \! \!  (1,3)_2\oplus(1,2)_2\oplus(1,1)_2\oplus(1,0)_1\oplus(0,4)_1\oplus(0,3)_1\oplus(0,2)_1\oplus(0,1)_1\oplus(2,3)_2\qquad \\
\nonumber \! \! \! \! \! &&\! \! \! \! \! \oplus(2,2)_2\oplus(2,1)_2\oplus(2,0)_1\oplus(1,4)_1\oplus(3,3)_1\oplus(3,2)_1\oplus(3,1)_1\oplus(3,0)_1\oplus(2,4)_1\eea
in terms of Dynkin components, and with the multiplicity appended as a subscript. 
}
\medskip
 
\subsection{Stretching}
The relation (\ref{K-br}) is also well suited for the study of the behaviour of branching coefficients under
``stretching". 
From (\ref{boundbr}) we learn that the growth is at most linear
$$ \mult_{s\alpha}(s{\bbeta}) \le s\min(\alpha_i-\alpha_{i+1}) +1\,.  $$

For example, for $n=3$, with Dynkin components,
$\mult_{(s,s)}(s)=s+1 $ since
$$   \mult_{(s,s)}(s)= \CKp( \{2s+2, s+1,0\}; \{s+1,0\}) = (s+1) \CKp(\{2,1,0\};\{1,0\})= s+1 $$
 while for $ \mult_{(s,s)}(s-1)$ or $ \mult_{(s,s)}(2s)$, we are not probing the function on its plateau and 
its behaviour is not always linear in $s$:
 $$\CKp (\{2s+2, s+1,0\}; \{{\bbeta}_1,0\})=\begin{cases}   {\bbeta}_1 & \mathrm{if}\quad 0\le {\bbeta}_1\le s+1\\
 2(s+1)-{\bbeta}_1 & \mathrm{if}\quad s+1\le {\bbeta}_1\le 2(s+1)  \end{cases}$$
 whence 
 $   \mult_{(s,s)}(s-1)= \CKp( \{2s+2, s+1,0\}; \{s,0\}) = s $ and 
 $ \mult_{(s,s)}(2s)= \CKp( \{2s+2, s+1,0\}; \{2s+1,0\}) = 1 $. 

 \begin{figure}[htb]
  \centering
       \includegraphics[width=10pc]{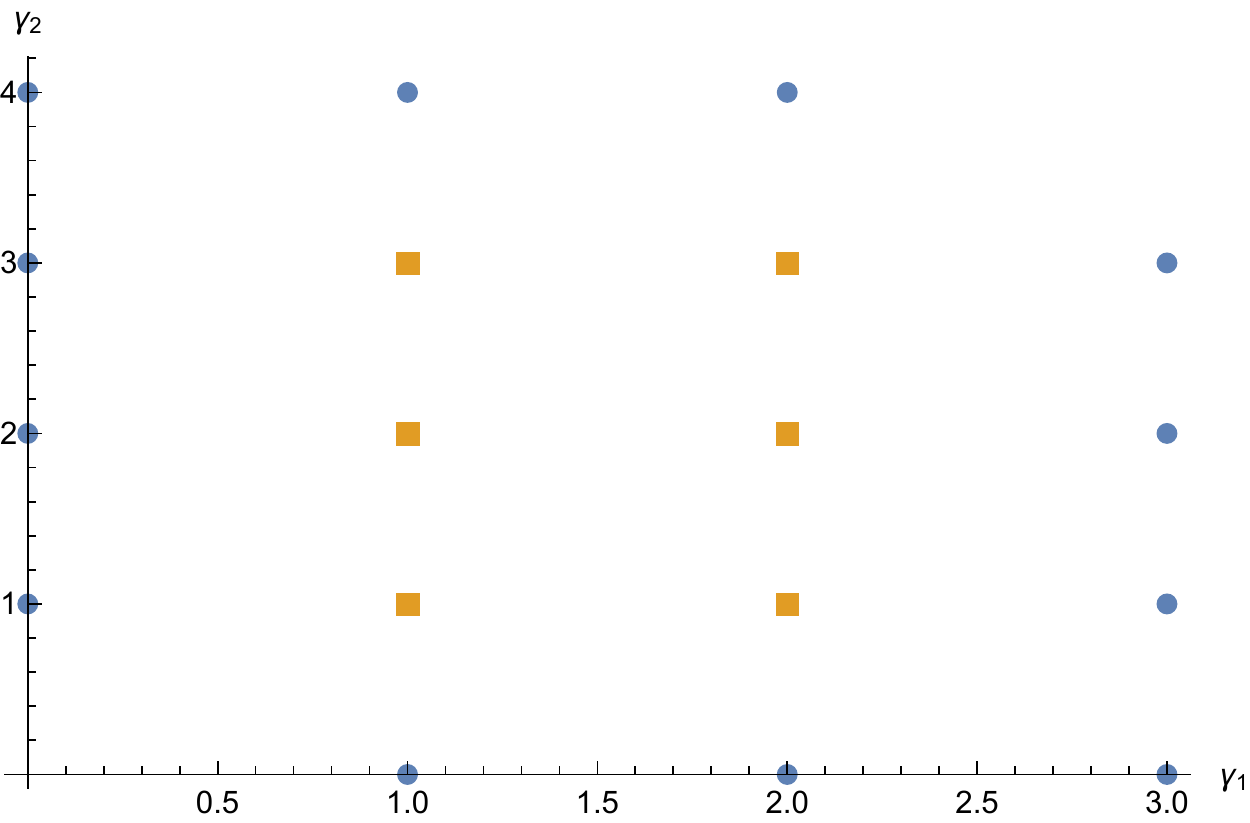}\qquad \includegraphics[width=10pc]{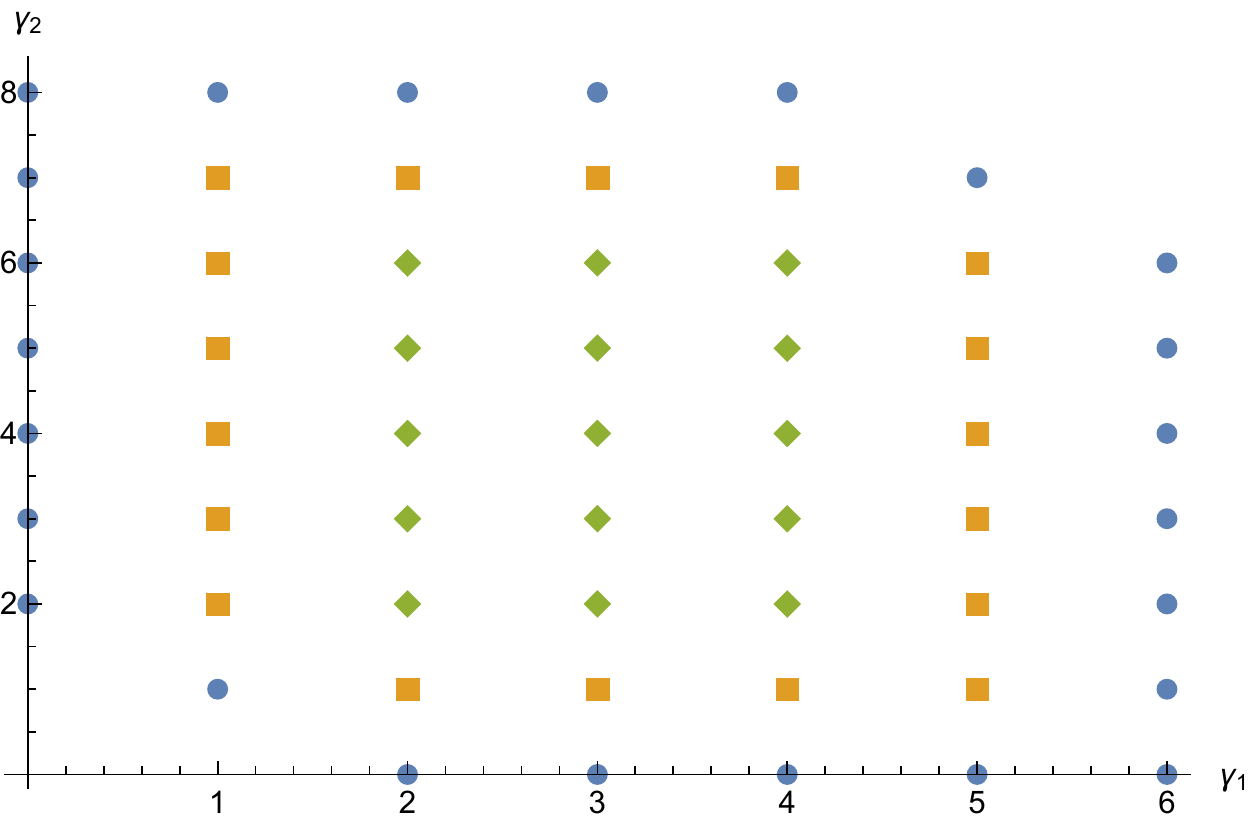}\qquad \includegraphics[width=12pc]{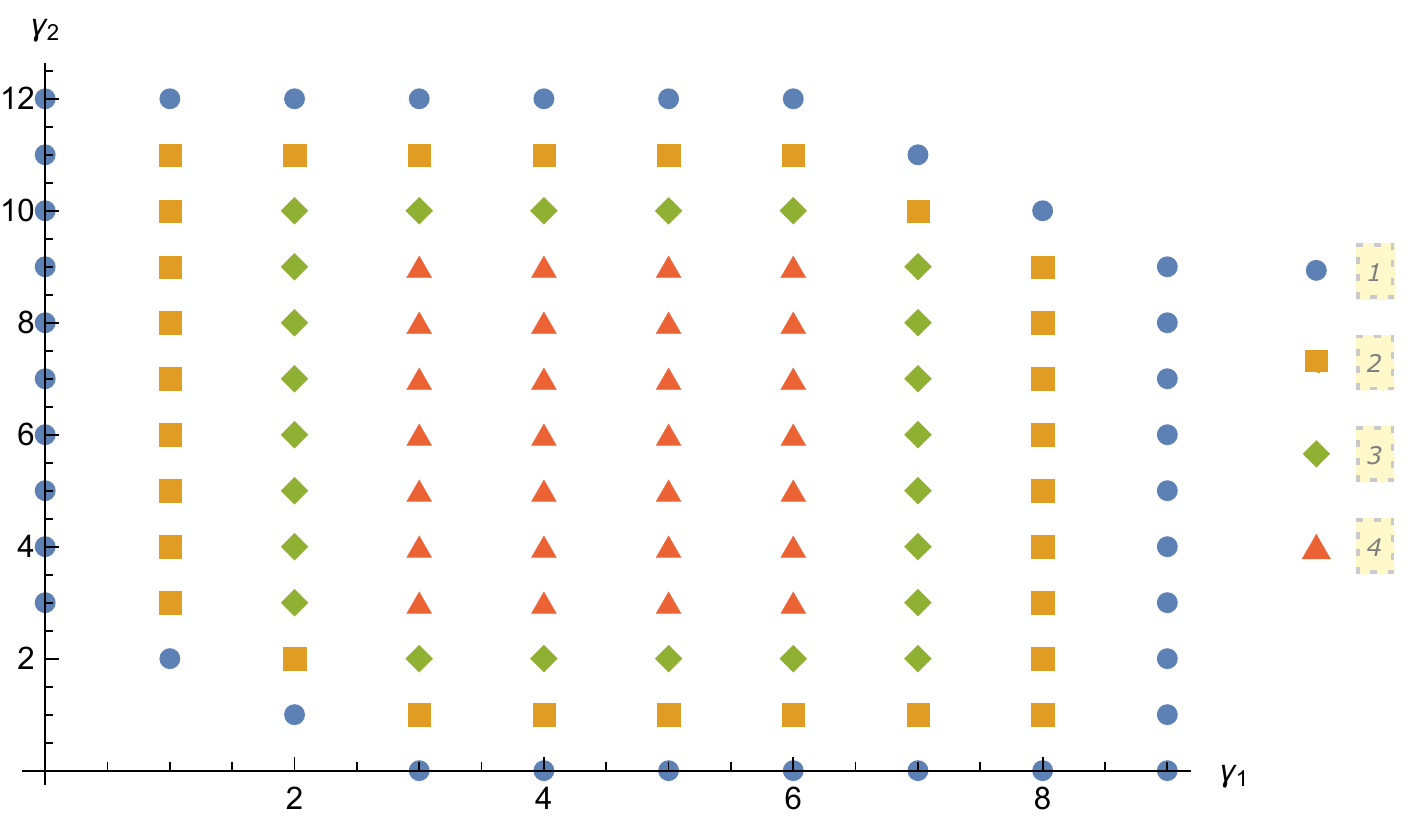}
       \caption{\label{multpatt} Weights in the ${\bbeta}$-plane (in Dynkin components) appearing in the decomposition of the weight $\alpha=s\{6,4,3,0\}\equiv 
       s(2,1,3)$ of $\SU(4)$, for $s=1,2,3$. Markers of different colours code for multiplicities from 1 to 4.}  
       \end{figure}

         \begin{figure}[htb]
  \centering
      \includegraphics[width=20pc]{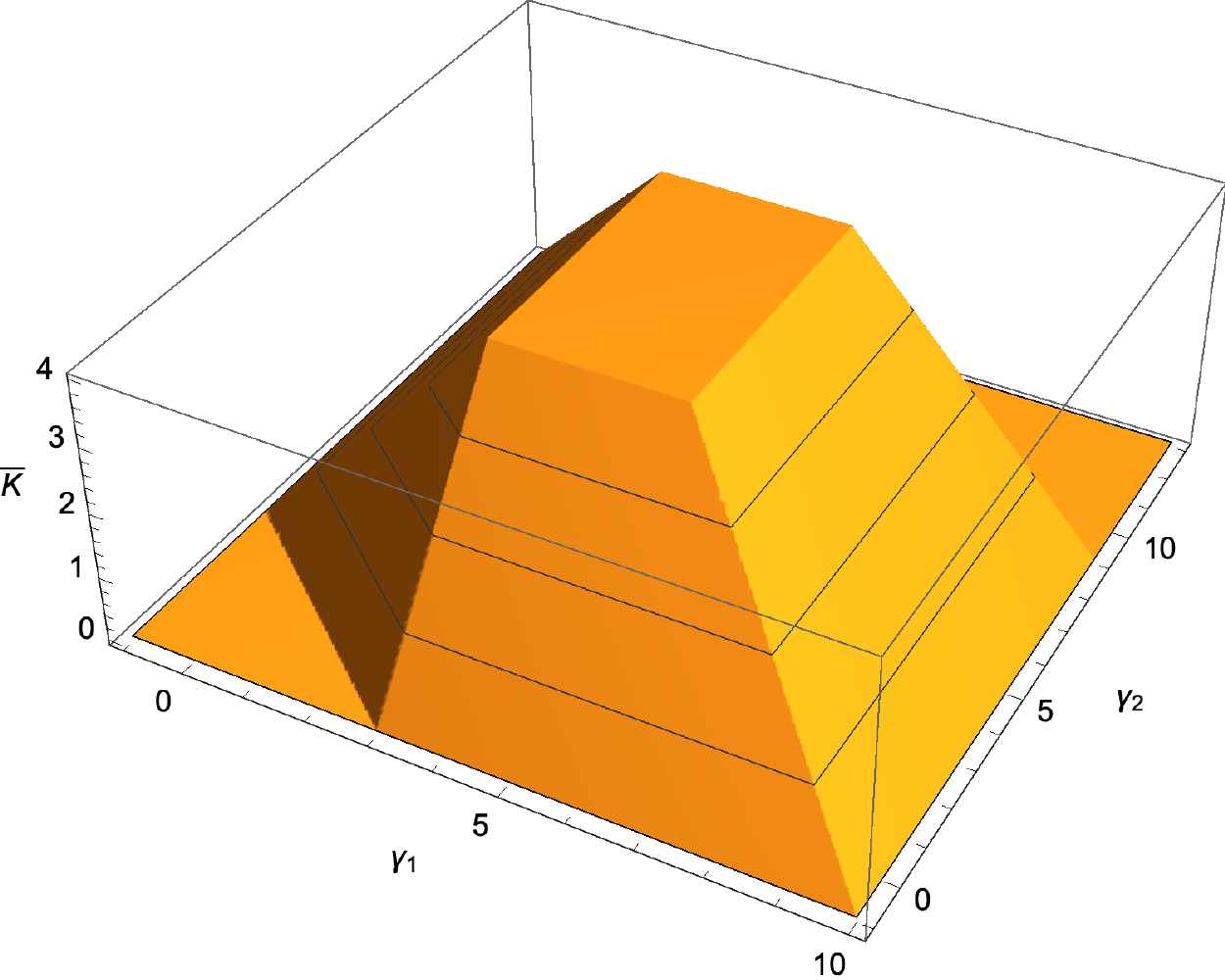}
       \caption{\label{figKp4b}  
       The $\CKp$ function  for $n=4$,   
             $\alpha+\rho_4=\{21,14,10,0\}=3\{6,4,3,0\}+\{3,2,1,0\}$ and ${\bbeta}_1,{\bbeta}_2$ are Dynkin components.
     The cross-sections at altitude $1,2,3,4$ match the successive layers of multiplicities in Fig. \ref{multpatt}, right. }   
       \end{figure}
    \newpage   
 
Similar behaviours occur for branching coefficients in higher rank cases, 
due to the linear growth of the maximal value (\ref{maxCKp}).
For $\SU(3)\subset\SU(4)$, the  points of increasing multiplicity form a matriochka pattern, see Fig. \ref{multpatt}, in a way already encountered 
in the Littlewood--Richardson coefficients of $\SU(3)$, \cite{CZ2014}.
This pattern just reproduces the cross-sections of increasing altitude of the Aztec pyramid
 of Fig. \ref{figKp4b}.


\section*{Acknowledgements} 
All my gratitude to Robert Coquereaux for his constant interest, encouragement and assistance. I also want to thank Jacques Faraut and Colin McSwiggen for their critical reading and comments.

  \end{document}